# Comparison of probabilistic and exact methods for estimating the asymptotic behavior of summation arithmetic functions

VICTOR VOLFSON

ABSTRACT. The paper compares probabilistic and exact methods for estimating the asymptotic behavior of summation arithmetic functions, and estimates of the results are obtained by precise methods. Conditions for stationarity in the broad sense are investigated for summation arithmetic functions. A lemma and theorems about the estimation of the standard deviation for the summation arithmetic Mertens $M(n)$ and Lowville $L(n)$ functions completely satisfying the stationarity conditions in the broad sense are proved.

## 1. INTRODUCTION

Based on Chebyshev's inequality, following [1], we can write the relation:

$$\lim_{n \to \infty} P(|f(n) - M[f(n)]| \leq h_1(n)\sqrt{D[f(n)]}) = 1, \qquad (1.1)$$

where $M[f(n)], D[f(n)]$ respectively, the mean value and variance of the arithmetic function $f(n)$, and $h_1(n)$ is an unboundedly increasing function $h_1(n) \leq n^\xi$, i.e. $h_1(n)$ is slowly growing function.

It follows from (1.1) that knowing the deviation from the mean value of an arithmetic function $D[f(n)]$ is very important.

We will understand a summation function as a sum of arithmetic functions:

$$S(x) = \sum_{n \leq x} f(n), \qquad (1.2)$$

_____________________________________________________________________ -

where $f(n)$ is the arithmetic function.

Examples of summation arithmetic functions are:

- the number of prime numbers not exceeding $x$ - $\pi(x) = \sum_{p \leq x} 1$, where $p$ - a prime number;

- Chebyshev functions - $\Psi(x) = \sum_{p^m < x} \log(p)$ and $\theta(x) = \sum_{p < x} \log(p)$;

- Mertens function - $M(x) = \sum_{n < x} \mu(n)$, where $\mu(n)$ is the Möbius function, $\mu(n) = 1$ if $n$ has an even number of prime factors of the first degree, $\mu(n) = -1$ if $n$ has an odd number of prime factors of the first degree, $\mu(n) = 0$ if $n$ has prime factors not only of the first degree;

- Liouville function - $L(x) = \sum_{n < x} \lambda(n)$, where $\lambda(n) = (-1)^{\nu(n)}$, where $\nu(n)$ is the number of prime factors in taking into account their multiplicity.

We will consider these summation arithmetic functions in the sequel.

An estimate is given for the variance for the summation function $S(n)$ in [2] under the stationary conditions in the broad sense:

$$D[S(n)] = nh_2(n), \qquad (1.3)$$

where $h_2(n)$ is a slowly growing function.

Having in mind (1.1) and (1.3), the following relation holds when the stationarity conditions are satisfied in the broad sense for $S(n)$, almost everywhere:

$$|S(n) - M[S(n)]| \leq h(n)n^{1/2}, \qquad (1.4)$$

where $h(n)$ is a slowly growing function.

It is shown in [3] that the summation arithmetic function will be stationary in the broad sense if the following conditions are satisfied:

$$1. \ S(n) = \sum_{k=1}^{n} f(k) = mn + o(n), \qquad (1.5)$$

where $m$ is a constant.

2. $|f(k)| \leq f$, (1.6)

where $f$ is a constant.

3. $\lim_{n \to \infty} \text{cov}[f(k)f(k+n)] = \lim_{n \to \infty} (M[f(k)f(k+n)] - M[f(k)]M[f(k+n)]) = 0$ (1.7)

Condition (1.7) corresponds to the asymptotic independence of the values of the arithmetic function $f(k), f(k+n)$.

We shall study the summation arithmetic functions in the paper, having the form (1.5). Summation arithmetic functions $\pi(x), M(x), L(x)$ satisfy condition (1.6). The condition (1.7) deserves to be explained, which is what we are going to do.

The summation arithmetic functions can be represented [3] as sums of random variables:
$$S_n = \sum_{i=1}^{n} x_i.$$ (1.8)

If the random variables $x_k$ and $x_{k+n}$ are asymptotically independent (at a value $n \to \infty$), then it is satisfied for any values $B_1, B_2 \in B$:

$$P(x_k \in B_1, x_{k+n} \in B_2) \to P(x_0 \in B_1)P(x_0 \in B_2).$$ (1.9)

Let's consider the arithmetic function of the number of prime numbers not exceeding the value $x$ - $\pi(x)$. Suppose, if the natural $k(k > 2)$ is prime, then the value of the random variable is $x_k = 1$, otherwise - $x_k = 0$. Similarly with respect to $x_{k+1}$. We denote: $P(x_k = 1) = a(0 \leq a \leq 1)$, $P(x_k = 0) = 1 - a$, $P(x_{k+1} = 1) = b(0 \leq b \leq 1)$,. Then, on the basis of (1.9), on the one hand:

$$P(x_k = 1, x_{k+1} = 0) = P(x_k = 1)P(x_{k+1} = 0 / x_k = 1) = a \cdot 1 = a,$$

because if $x_k = 1$, then $k$ is an odd number. On the other hand: $P(x_k = 1)P(x_{k+1} = 0) = a(1-b)$. Therefore: $P(x_k = 1, x_{k+1} = 0) = a$ is not equal to $P(x_k = 1)P(x_{k+1} = 0) = a(1-b)$, if $b$ not equal to $0$.

Thus, random variables $x_k, x_{k+1}$ are dependent. This situation persists with a value $n \to \infty$. Consequently, the asymptotic independence condition does not hold for the summation function $\pi(x)$.

Asymptotic independence is satisfied for summation arithmetic functions, the defining property of which is the presence of an even or odd number of prime divisors for a natural number, with a large number of prime divisors, i.e. with the value $n \to \infty$. This applies to the arithmetic summation functions $M(x), L(x)$.

We will try to prove the results obtained above with the help of probabilistic methods for estimating the asymptotic behavior of summation arithmetic functions in the following sections using exact methods.

## 2. INVESTIGATION OF THE ASYMPTOTIC BEHAVIOR OF SUMMATIONARITHMETIC FUNCTIONS USING EXACT METHODS

We will carry out a study of the asymptotic behavior of the summation arithmetic functions: $M(n) = \sum_{k=1}^{n} \mu(k), L(n) = \sum_{k=1}^{n} \lambda(k)$ in this chapter. We begin with the proof of the lemma for the arithmetic functions of Liouville $\lambda(k)$ and Mobius $\mu(k)$.

Lemma 1. The following estimate holds for the Liouville and Mobius arithmetic functions:

$$\sum_{i=1}^{n} \sum_{j=1}^{n} f(i) f(j) = o(n^2).  \qquad (2.1)$$

Proof

It is proved on page 87 [4] that there is an asymptotically identical number of natural numbers with an even and odd number of prime factors.

We denote these: $N_1(x) = \sum_{n \leq x, \lambda(n)=1} 1, N_{-1}(x) = \sum_{n \leq x, \lambda(n)=-1} 1$.

This can be written in the form:

$$N_1(x) = x/2 + o(x), N_{-1}(x) = x/2 + o(x).  \qquad (2.2)$$

The product is equal to $\lambda(i)\lambda(j) = 1$, if $\lambda(i) = 1, \lambda(j) = 1$ or $\lambda(i) = -1, \lambda(j) = -1$.

We denote $N_{11}(x) = \sum_{i,j \leq x, \lambda(i)=1, \lambda(j)=1} 1, N_{-1-1}(x) = \sum_{i,j \leq x, \lambda(i)=-1, \lambda(j)=-1} 1$.

Having in mind (2.2) it can be written:

$N_{11} = x^2/4 + o(x^2)$, $N_{-1-1} = x^2/4 + o(x^2)$, thus

$$N_{11} + N_{-1-1} = x^2/2 + o(x^2). \qquad (2.3)$$

Similarly to (2.3), we get the number of natural numbers for which $\lambda(i)\lambda(j) = -1$:

$$N_{1-1} + N_{-11} = x^2/2 + o(x^2). \qquad (2.4)$$

On the basis of (2.3) and (2.4), the formula (2.1) holds for the arithmetic Liouville function.

A similar proof is valid for the Mobius function, taking into account that for natural numbers there are not squares-free values $\mu(n) = 0$ and therefore the corresponding products are zero. Consequently, these products do not affect the amount of the sum $\sum_{i=1}^{n}\sum_{j=1}^{n}\mu(i)\mu(j)$.

There is an asymptotic equal number of natural numbers having an even and odd numbers of prime factors among the natural numbers that are free from squares, so the proof is similar to that for the Lowville arithmetic function.

Theorem 1

An estimate is performed for the Liouville and Mobius arithmetic functions:

$$M[f(i)f(j)] - M[f(i)]M[f(j)] = O(1/n), \qquad (2.5)$$

where $M[...]$ is the average value of the arithmetic function enclosed in parentheses.

Proof

$$M[f(i)f(j)] = \frac{\sum\sum_{i \neq j} f(i)f(j)}{n(n-1)} < \frac{\sum_{i=1}^{n}\sum_{j=1}^{n} f(i)f(j)}{n(n-1)}.$$

We will take into account that:

$$\frac{\sum_{i=1}^{n}\sum_{j=1}^{n} f(i)f(j)}{n(n-1)} = \frac{(\sum_{k=1}^{n} f(k))^2}{n(n-1)} = \frac{\sum\sum_{i \neq j} f(i)f(j)}{n(n-1)} + O(1/n) \qquad (2.6)$$

$$M[f(i)]M[f(j)] = \frac{\sum_{i=1}^{n} f(i)}{n} \cdot \frac{\sum_{j=1}^{n} f(i)}{n} = \frac{(\sum_{i=1}^{n} f(i))^2}{n^2}. \qquad (2.7)$$

Based on (2.6) and (2.7) we obtain:

$$M[f(i)f(j)] - M[f(i)]M[f(j)] = \frac{(\sum_{k=1}^{n} f(k))^2}{n(n-1)} - \frac{(\sum_{i=1}^{n} f(k))^2}{n^2} + O(1/n) = \frac{(\sum_{k=1}^{n} f(k))^2}{n^2(n-1)} + O(1/n).$$

Therefore, based on Lemma 1 we obtain:

$$M[f(i)f(j)] - M[f(i)]M[f(j)] = \frac{(\sum_{k=1}^{n} f(k))^2}{n^2(n-1)} + O(1/n) = \frac{o(n^2)}{n^2(n-1)} + O(1/n) = O(1/n).$$

Theorem 1 will be used in the proof of the following Theorem 2.

The estimate (2.5) has an independent value, since it proves the asymptotic independence (1.7) for the arithmetic functions of Liouville and Mobius. Thus, these functions satisfy all the stationarity conditions in the broad sense (1.5), (1.6), (1.7).

Theorem 2

The order of the standard deviation of the summation arithmetic functions $L(n), M(n)$ is $O(n^{1/2})$.

Proof

Let's consider the deviation of the summation arithmetic function $S(n)$ from its mean value:

$$F(n) = \sum_{k=1}^{n} f(k) - M[\sum_{k=1}^{n} f(k)] = \sum_{k=1}^{n} (f(k) - M[f(k)]) = \sum_{k=1}^{n} a(k), \qquad (2.8)$$

where $a(k) = f(k) - M[f(k)]$.

Having jn mind (2.8) we define the arithmetic function $F^2(n)$:

$$F^2(n) = (\sum_{k=1}^{n} a(k))^2 = \sum_{k=1}^{n} a^2(k) + \sum_{i \neq j} a(i)a(j). \qquad (2.9)$$

Taking (2.9) into account, we find the mean value of the arithmetic function $F^2(n)$:

$$M[F^2(n)] = M[\sum_{k=1}^{n} a^2(k)] + M[\sum_{i \neq j} a(i)a(j)] \qquad (2.10)$$

We estimate the first term in (2.10):

$$M[\sum_{k=1}^{n} a^2(k)] = \sum_{k=1}^{n} M[(f(k) - M[f(k)])^2] = \sum_{k=1}^{n} (M[f^2(k)] - M^2[f(k)]) \leq \sum_{k=1}^{n} M[f^2(k)].$$

Taking into account that for arithmetic functions $\lambda(k), \mu(k)$ is satisfied $|f(k)|\leq 1$ and $M[f^2(k)]\leq 1$ we get:

$$M[\sum_{k=1}^{n} a^2(k)] \leq \sum_{k=1}^{n} M[f^2(k)] \leq n, \qquad (2.11)$$

those $M[\sum_{k=1}^{n} a^2(k)] = O(n)$.

We estimate the second term in (2.10):

$$M[\sum_{i\neq j} a(i)a(j)] = M\sum_{i\neq j}(f(i)-M[f(i)])(f(j)-M[f(j)]) = \sum_{i\neq j}(M[f(i)f(j)]-M[f(i)]M[f(j)])$$

Based on Theorem 1:

$$M[\sum_{i\neq j} a(i)a(j)] = n(n-1)O(1/n) = O(n). \qquad (2.12)$$

Having in mind (2.11) and (2.12), the order of the mean value of the arithmetic function $F^2(n)$ is:

$$M[F^2(n)] = O(n). \qquad (2.13)$$

We obtain (from (2.13)) the required estimate for the order of the standard deviation for the summation arithmetic functions $L(n), M(n)$:

$$\sigma[S(n)] = O(n^{1/2}). \qquad (2.14)$$

An estimate similar to (2.14) was obtained in [5] with respect to other summation arithmetic functions (statistical mechanics). In the same paper, the concept of asymptotic independence is used, as the main property of the summands of arithmetic functions, which gives such an estimate. This gives me confidence in the result.

3. CONCLUSION AND SUGGESTIONS FOR FURTHER WORK

Let's draw the following conclusions:

1. We prove the estimate of the standard deviation for the summation arithmetic functions of Liouville and Mertens (Theorem 2) using the exact methods:

    $$\sigma[S(n)] = O(n^{1/2}). \qquad (3.1)$$

2. Based on (1.1) using (3.1), the estimate for the deviation of the summation arithmetic functions of Liouville and Mertens from their average value is valid almost everywhere:

$$F(n) = O(n^{1/2} \varphi(n)),  \quad (3.2)$$

where $\varphi(n)$ is any slowly increasing function.

3. Formula (3.2) confirms the estimate (1.4) of the work for the summation arithmetic functions of Liouville and Mertens that satisfy the stationarity conditions in the broad sense (when conditions (1.5), (1.6), (1.7) are satisfied).

4. The condition of asymptotic independence (1.7) for the indicated arithmetic functions was proved with the help of exact methods (Theorem 1).

The next article will continue to study the asymptotic behavior of arithmetic functions.

4. ACKNOWLEDGEMENTS

Thanks to everyone who has contributed to the discussion of this paper. I am grateful to everyone who expressed their suggestions and comments in the course of this work.